\documentclass[11pt,a4paper]{article}
\usepackage{pdfsync} 

\usepackage{amsfonts,amsmath,amsthm}
\newtheorem{theorem}{Theorem}
\numberwithin{theorem}{section}
\newtheorem{lemma}[theorem]{Lemma}
\newtheorem{corollary}[theorem]{Corollary}
\newtheorem{example}[theorem]{Example}
\numberwithin{equation}{section}
\usepackage{amssymb}
\usepackage{hyperref}
\usepackage{graphicx}
\usepackage{wrapfig}
\usepackage{subfig}
\usepackage{epsfig}
\usepackage{float,mathrsfs}
\usepackage{bbm}

\usepackage{algpseudocode}
\usepackage[font=scriptsize]{caption}
\usepackage{chngcntr} 
\counterwithout{equation}{section}
\counterwithout{figure}{section}
\usepackage{color}
\definecolor{refkey}{rgb}{0.9451,0.2706,0.4941}\definecolor{labelkey}{rgb}{0.9451,0.2706,0.4941}

\definecolor{darkred}{RGB}{139,0,0}
\definecolor{darkgreen}{RGB}{0,100,0}
\definecolor{darkmagenta}{RGB}{139,0,139}

\newcommand{\mask}[1]{{}}

\title{Bounds on the spectrum of nonsingular triangular $(0,1)$-matrices}

\author{V. Kaarnioja\footnotemark[2]}

\renewcommand{\thefootnote}{\fnsymbol{footnote}}

\begin{document}
\maketitle

\begin{abstract}
Let $K_n$ be the set of all nonsingular $n\times n$ lower triangular $(0,1)$-matrices. Hong and Loewy (2004) introduced the numbers
$$
c_n={\rm min}\{\lambda\mid \lambda~\text{is an eigenvalue of}~XX^{\rm T},~X\in K_n\},\quad n\in\mathbb{Z}_+.
$$
A related family of numbers was considered by Ilmonen, Haukkanen, and Merikoski (2008):
$$
C_n={\rm max}\{\lambda\mid \lambda~\text{is an eigenvalue of}~XX^{\rm T},~X\in K_n\},\quad n\in\mathbb{Z}_+.
$$
These numbers can be used to bound the singular values of matrices belonging to $K_n$ and they appear, e.g., in eigenvalue bounds for power GCD matrices, lattice-theoretic meet and join matrices, and related number-theoretic matrices. In this paper, it is shown that for $n$ odd, one has the lower bound
$$
c_n\geq \frac{1}{\sqrt{\frac{1}{25}\varphi^{-4n}+\frac{2}{25}\varphi^{-2n}-\frac{2}{5\sqrt{5}}n\varphi^{-2n}-\frac{23}{25}+n+\frac{2}{25}\varphi^{2n}+\frac{2}{5\sqrt{5}}n\varphi^{2n}+\frac{1}{25}\varphi^{4n}}},
$$
and for $n$ even, one has
$$
c_n\geq \frac{1}{\sqrt{\frac{1}{25}\varphi^{-4n}+\frac{4}{25}\varphi^{-2n}-\frac{2}{5\sqrt{5}}n\varphi^{-2n}-\frac{2}{5}+n+\frac{4}{25}\varphi^{2n}+\frac{2}{5\sqrt{5}}n\varphi^{2n}+\frac{1}{25}\varphi^{4n}}},
$$
where $\varphi$ denotes the golden ratio. These lower bounds improve the estimates derived previously by Mattila (2015) and Alt{\i}n{\i}\c{s}{\i}k et al.~(2016). The sharpness of these lower bounds is assessed numerically and it is conjectured that $c_n\sim 5\varphi^{-2n}$ as $n\to\infty$. In addition, a new closed form expression is derived for the numbers $C_n$, viz.
$$
C_n=\frac14 \csc^2\bigg(\frac{\pi}{4n+2}\bigg)=\frac{4n^2}{\pi^2}+\frac{4n}{\pi^2}+\bigg(\frac{1}{12}+\frac{1}{\pi^2}\bigg)+\mathcal{O}\bigg(\frac{1}{n^2}\bigg),\quad n\in\mathbb{Z}_+.
$$
\end{abstract}

\footnotetext[2]{School of Mathematics and Statistics, University of New South Wales, Sydney NSW 2052, Australia ({\tt vesa.kaarnioja@iki.fi}). Present address: School of Engineering Science, LUT University, P.O. Box 20, FI-53851 Lappeenranta, Finland.}

\renewcommand{\thefootnote}{\arabic{footnote}}

\section{Introduction}

Let $K_n$ denote the set of all nonsingular $n\times n$ lower triangular $(0,1)$-matrices. For example, $K_3$ consists of the elements
\begin{align*}
&\begin{pmatrix}
1&0&0\\
0&1&0\\
0&0&1
\end{pmatrix},~\begin{pmatrix}
1&0&0\\
1&1&0\\
0&0&1
\end{pmatrix},~\begin{pmatrix}
1&0&0\\
0&1&0\\
1&0&1
\end{pmatrix},~\begin{pmatrix}
1&0&0\\
0&1&0\\
0&1&1
\end{pmatrix},\\
&\begin{pmatrix}
1&0&0\\
1&1&0\\
1&0&1
\end{pmatrix},~\begin{pmatrix}
1&0&0\\
1&1&0\\
0&1&1
\end{pmatrix},~\begin{pmatrix}
1&0&0\\
0&1&0\\
1&1&1
\end{pmatrix},~\begin{pmatrix}
1&0&0\\
1&1&0\\
1&1&1
\end{pmatrix},
\end{align*}
and it is easy to see that $\#K_n=2^{n(n-1)/2}$ for all $n\in\mathbb{Z}_+$.

Hong and Loewy~\cite{hongloewy} introduced the numbers
\begin{align*}
c_n=\min\{\lambda\mid \lambda~\text{is an eigenvalue of}~XX^{\rm T},~X\in K_n\},\quad n\in\mathbb{Z}_+,
\end{align*}
as a means to give a lower bound for the smallest eigenvalue of power GCD matrices defined on any set of positive integers. A closely related sequence of numbers 
\begin{align*}
C_n=\max\{\lambda\mid \lambda~\text{is an eigenvalue of}~XX^{\rm T},~X\in K_n\},\quad n\in\mathbb{Z}_+,
\end{align*}
was introduced by Ilmonen, Haukkanen, and Merikoski~\cite{ihm} in order to derive upper bounds for the largest eigenvalues of lattice-theoretic meet and join matrices.

The numbers $c_n$ and $C_n$ have an intimate connection with the extremal singular values of matrices belonging to $K_n$. Let $\sigma_{\min}(X)$ and $\sigma_{\max}(X)$ denote the smallest and largest singular values of matrix $X$, respectively. Then
$$
\min_{X\in K_n}\sigma_{\min}(X)=\sqrt{c_n}\quad\text{and}\quad \max_{X\in K_n}\sigma_{\max}(X)=\sqrt{C_n},\quad n\in\mathbb{Z}_+.
$$

The following example illustrates how the numbers $c_n$ and $C_n$ can be used to give bounds for the smallest and largest eigenvalues of certain number-theoretic matrices.
\begin{example}[{cf.~\cite{ihm}}]\rm
Let $(P,\preceq,\wedge,\hat 0)$ be a locally finite meet semilattice, where $\preceq$ is a partial ordering on the set $P$, $\wedge$ denotes the {\em meet} (or {\em greatest lower bound}) of two elements in $P$, and $\hat 0\in P$ is the least element such that $\hat 0\preceq x$ for all $x\in P$. Let $S=\{x_1,\ldots,x_n\}\subset P$ be a lower closed set such that $x_i\preceq x_j$ only if $i\leq j$. Let $f\!:P\to\mathbb{R}$ be a function and define the $n\times n$ meet matrix $A$ elementwise by setting $A_{i,j}=f(x_i\wedge x_j)$ for $i,j\in\{1,\ldots,n\}$. Define the function $$J_{P,f}(x)=\sum_{\hat 0\preceq z\preceq x}f(z)\mu(z,x)\quad\text{for all}~x\in P,$$ where $\mu$ denotes the M\"{o}bius function of $P$. If $J_{P,f}(x)>0$ for all $x\in S$, then
\begin{align}
\lambda_{\min}(A)\geq c_n\min_{x\in S}J_{P,f}(x)\quad\text{and}\quad \lambda_{\max}(A)\leq C_n\max_{x\in S}J_{P,f}(x).\label{eq:application}
\end{align}
For example, in the case of the divisor lattice $(\mathbb{Z}_+,|,{\rm gcd})$, $f(x)=x^\alpha$ for $x\in\mathbb{Z}_+$, and $\alpha>0$,  the function $J_{P,f}$ coincides with the Jordan totient function and the matrix $A$ is called a power GCD matrix. (If $\alpha=1$, then $J_{P,f}$~is Euler's totient function and the matrix $A$ is called a GCD matrix.) See~\cite{hongloewy,ihm} for further discussion and generalizations of this result.
\end{example}

The eigenvalues of power GCD and LCM matrices, meet and join matrices, and a variety of closely related number-theoretic matrices have been considered by many authors. Balatoni~\cite{balatoni} gave lower and upper bounds on the extremal eigenvalues of Smith's matrix, i.e., the matrix with its $(i,j)$ element equal to ${\rm gcd}(i,j)$, which was considered in 1875/76 by the eponymous Smith~\cite{smith}. Beslin and Ligh showed that GCD matrices defined on an arbitrary set of distinct positive integers are positive definite~\cite{beslinligh}. Bourque and Ligh investigated GCD matrices composed elementwise with arithmetic functions and determined conditions which ensure the positive definiteness of these matrices~\cite{bourqueligh}. Further progress on this line of research was made by Hong~\cite{hong2008}, who studied the behavior of the largest eigenvalue of GCD matrices associated with certain arithmetic functions, while research on the behavior of the corresponding smallest eigenvalue was conducted by Hong and Loewy~\cite{hongloewy11}. Hong and Lee obtained results on the asymptotic behavior of the eigenvalues of reciprocal power LCM matrices~\cite{hongenochlee2008}. Recently, Merikoski revisited the lower bound of the smallest eigenvalue of Smith's matrix~\cite{merikoski}, Mattila and Haukkanen derived eigenvalue bounds for ``mixed'' power GCD and LCM matrices~\cite{haukkanenmattila12}, Alt{\i}n{\i}\c{s}{\i}k and B\"{u}y\"{u}kk\"{o}se derived new eigenvalue bounds for GCD and LCM matrices~\cite{altinisik16}, and Ilmonen obtained bounds for the eigenvalues of meet hypermatrices~\cite{ilmonen}. Haukkanen et al.~gave an interesting lattice-theoretic interpretation for the inertia of LCM matrices~\cite{inertia}.

In their work~\cite{hongloewy}, Hong and Loewy did not give a lower bound or other estimates for $c_n$ and there do not appear to have been any developments on estimating $c_n$~in the literature until Alt{\i}n{\i}\c{s}{\i}k and B\"{u}y\"{u}kk\"{o}se~\cite{altinisik15} analyzed a related quantity, which could be used to obtain an upper bound on the numbers $c_n$. Mattila~\cite{mattila} derived the following lower bounds for $c_n$:
\begin{align}
c_n&\geq \bigg(\frac{48}{n^4+56n^2+48n}\bigg)^{(n-1)/2}\quad\text{for even}~n,\label{eq:mattila1}\\
c_n&\geq \bigg(\frac{48}{n^4+50n^2+48n-51}\bigg)^{(n-1)/2}\quad\text{for odd}~n.\label{eq:mattila2}
\end{align}
The lower bounds~\eqref{eq:mattila1} and~\eqref{eq:mattila2} were subsequently improved in~\cite{altinisik}:
\begin{align}
c_n\geq \frac{2}{2F_nF_{n+1}+(-1)^n+1}\quad\text{for}~n\in\mathbb{Z}_+,\label{eq:altinisik}
\end{align}
where $(F_n)_{n=1}^\infty$ denotes the Fibonacci sequence. However, a straightforward numerical investigation shows that the bounds~\eqref{eq:mattila1}--\eqref{eq:altinisik} are {\em not} sharp. It is the goal of this article  to remedy this situation by developing an improved lower bound for the numbers $c_n$. In addition, a new characterization for the numbers $C_n$ is also derived in this paper.

This paper is structured as follows. Section~\ref{sec:hongloewy} begins with the development of a new lower bound on Hong and Loewy's numbers $c_n$. In Subsection~\ref{sec:main}, it is shown that this lower bound can be expressed in a much simplified form, which is the main contribution of this paper. The sharpness of this new lower bound is assessed by numerical experiments in Subsection~\ref{sec:numex} along with numerical comparisons involving the bounds~\eqref{eq:mattila1}--\eqref{eq:altinisik}. A novel characterization for the closely related sequence of  Ilmonen--Haukkanen--Merikoski numbers $C_n$ is proved in Section~\ref{sec:ihm}. Finally, some conclusions and thoughts about future work are given at the end of the paper.

\section{Hong and Loewy's numbers $c_n$}\label{sec:hongloewy}

Alt{\i}n{\i}\c{s}{\i}k et al.~\cite{altinisik} proved the following characterization
\begin{align}
c_n=\lambda_{\min}(Z_n^{-1})\quad\text{for all }n\in\mathbb{Z}_+,\label{eq:ihm}
\end{align}
where $Z_n$ is the symmetric $n\times n$ matrix defined elementwise by
$$
(Z_n)_{i,j}=\begin{cases}
1+\sum_{k=i+1}^nF_{k-i}^2&\text{if }i=j,\\
(-1)^{j-i}(F_{j-i}+\sum_{k=j+1}^nF_{k-i}F_{k-j})&\text{if }i<j,\\
(-1)^{i-j}(F_{i-j}+\sum_{k=i+1}^nF_{k-i}F_{k-j})&\text{if }i>j,
\end{cases}
$$
for $i,j\in\{1,\ldots,n\}$ and the sequence of \emph{Fibonacci numbers} is defined by the recurrence relation $F_0=0$, $F_1=1$, and $F_k=F_{k-1}+F_{k-2}$ for $k\geq 2$.

The following technical result  will serve as the basis for the analysis in Subsection~\ref{sec:simplifiedbound}.
\begin{lemma}\label{lem:technicallemma} It holds for all $n\in\mathbb{Z}_+$ that
$$
c_n\geq \frac{1}{\sqrt{1+\sum_{i=2}^n(1+F_iF_{i-1})^2+2\sum_{i=2}^n\sum_{j=2}^i\big(F_{j-1}+\sum_{k=j+1}^iF_{k-1}F_{k-j}\big)^2}}.
$$
\end{lemma}
\proof Let $n\in\mathbb{Z}_+$. By the characterization~\eqref{eq:ihm}, it holds that
\begin{align}
c_n=\lambda_{\max}(Z_n)^{-1}=\|Z_n\|_2^{-1}\geq \|Z_n\|_{\rm F}^{-1},\label{eq:ZnFrob}
\end{align}
where $\|\cdot\|_2$ denotes the spectral norm, $\|\cdot\|_{\rm F}$ is the  Frobenius norm, and the final inequality is due to  $\|\cdot\|_2\leq \|\cdot\|_{\rm F}$. To prove the claim, it is sufficient to compute the value of the Frobenius norm appearing in~\eqref{eq:ZnFrob}.

Making use of the block structure of the matrices $Z_n$, it is possible to write
$$
Z_1=(1)\quad\text{and}\quad Z_n=\begin{pmatrix}a_n&b_n^{\rm T}\\
b_n&Z_{n-1}\end{pmatrix}\quad \text{for}~n\geq 1,
$$
where $a_n=(Z_n)_{1,1}=1+F_nF_{n-1}$ and the $(n-1)$-vector $b_n=[(Z_n)_{2,1},\ldots,(Z_n)_{n,1}]^{\rm T}\in\mathbb{R}^{n-1}$ clearly satisfies
$$
\|b_n\|^2=\sum_{j=2}^n\bigg(F_{j-1}+\sum_{k=j+1}^nF_{k-1}F_{k-j}\bigg)^2.
$$
Hence $\|Z_n\|_{\rm F}^2=a_n^2+2\|b_n\|^2+\|Z_{n-1}\|_{\rm F}^2$, which yields the recurrence relation
\begin{align*}
&\|Z_1\|_{\rm F}^2=1,\\
&\|Z_n\|_{\rm F}^2=\|Z_{n-1}\|_{\rm F}^2+(1+F_nF_{n-1})^2+2\sum_{j=2}^n\bigg(F_{j-1}+\sum_{k=j+1}^nF_{k-1}F_{k-j}\bigg)^2,\quad n\geq 2.
\end{align*}
This recurrence can be used to produce the expression
\begin{align*}
\|Z_n\|_{\rm F}^2=1+\sum_{i=2}^n(1+F_iF_{i-1})^2+2\sum_{i=2}^n\sum_{j=2}^i\bigg(F_{j-1}+\sum_{k=j+1}^iF_{k-1}F_{k-j}\bigg)^2,\quad n\in\mathbb{Z}_+,
\end{align*}
which, together with the inequality~\eqref{eq:ZnFrob}, proves the assertion.\endproof

Lemma~\ref{lem:technicallemma} gives a computable, albeit rather unwieldy, lower bound on the numbers $c_n$. However, it is shown in the following section that this lower bound can be recast into a much simpler closed form expression.
\subsection{Simplifying the lower bound on $c_n$}\label{sec:main}\label{sec:simplifiedbound}
In this section, a closed form expression for the term inside the square root in Lemma~\ref{lem:technicallemma} is derived. To this end, recall that the sequence of \emph{Lucas numbers} can be defined by the recursion $L_0=2$, $L_1=1$, and $L_k=L_{k-1}+L_{k-2}$ for $k\geq 2$. The Fibonacci--Binet formula and the Lucas--Binet formula can be used to write the Fibonacci and Lucas numbers explicitly as
\begin{align}
F_k=\frac{\varphi^k-(-\varphi)^{-k}}{\sqrt{5}}\quad\text{and}\quad L_k=\varphi^k+(-\varphi)^{-k}\quad\text{for all}~k\in\mathbb{Z}_{\geq 0},\label{eq:binet}
\end{align}
where $\varphi$ denotes the golden ratio.

The main result of this paper is given by the following theorem. It is a simplified version of the lower bound given in Lemma~\ref{lem:technicallemma}.

\begin{theorem}\label{thm:main}
It holds for all $n\in\mathbb{Z}_+$ that
$$
c_n\!\geq\! \frac{1}{\sqrt{\frac{1}{25}\varphi^{-4n}\!+\!\frac{3+(-1)^n}{25}\varphi^{-2n}\!-\!\frac{2}{5\sqrt{5}}n\varphi^{-2n}\!+\!\frac{13(-1)^n-33}{50}\!+\!n\!+\!\frac{3+(-1)^n}{25}\varphi^{2n}\!+\!\frac{2}{5\sqrt{5}}n\varphi^{2n}\!+\!\frac{1}{25}\varphi^{4n}}}.
$$
\end{theorem}
\proof
The proof is based on simplifying the term inside the square root in Lemma~\ref{lem:technicallemma}.

The claim is clearly true with equality for $n=1$. In the following analysis, let $n\geq i\geq j\geq 2$ be integers. Using the formulae~\eqref{eq:binet},  it is straightforward to check that
$$
F_{j-1}+\sum_{k=j+1}^iF_{k-1}F_{k-j}=\frac15\bigg(L_{2i-j}+\frac52 F_{j-1}-\frac12 (-1)^{i-j}L_{j-1}\bigg),
$$
where the sum is taken to be $0$ if the index set is empty. In consequence, it follows that
\begin{align}
\begin{array}{ll}
\bigg(F_{j-1}+\displaystyle\sum_{k=j+1}^iF_{k-1}F_{k-j}\bigg)^2\!\!\!&=\displaystyle\frac{1}{25}L_{2i-j}^2+\displaystyle\frac15 F_{j-1}L_{2i-j}-\displaystyle\frac{1}{25}(-1)^{i-j}L_{j-1}L_{2i-j}+\displaystyle\frac14 F_{j-1}^2\\
&\quad -\displaystyle\frac{1}{10}(-1)^{i-j}F_{j-1}L_{j-1}+\displaystyle\frac{1}{100}L_{j-1}^2.
\end{array}\label{eq:innerinterior}
\end{align}
Using the summation formulae
\begin{align*}
\begin{split}
&\sum_{j=2}^iL_{2i-j}^2=L_{2i-2}L_{2i-1}-L_{i}L_{i-1}\\
&\sum_{j=2}^iF_{j-1}L_{2i-j}=(i\!-\!1)F_{2i-1}\!+\!\frac{L_{2i-2}}{5}\!+\!\frac{2}{5}(-1)^i\\
&\sum_{j=2}^i(-1)^jL_{j-1}L_{2i-j}=2\!+\!\frac{1\!+\!(-1)^i}{2}L_{2i-1}\!-\!L_{2i-2}
\end{split}
\begin{split}
&\sum_{j=2}^iF_{j-1}^2=F_{i-1}F_i\\
&\sum_{j=2}^i(-1)^jF_{j-1}L_{j-1}=(-1)^i\frac{F_{i+1}^2-F_{i-2}^2}{4}\\
&\sum_{j=2}^iL_{j-1}^2=L_{i-1}L_i-2
\end{split}
\end{align*}
in conjunction with~\eqref{eq:innerinterior} yields that
\begin{align*}
\sum_{j=2}^i\bigg(F_{j-1}+\sum_{k=j+1}^iF_{k-1}F_{k-j}\bigg)^2=&\:\frac{1}{25}L_{2i-2}L_{2i-1}-\frac{3}{100}L_iL_{i-1}+\frac15 iF_{2i-1}-\frac15  F_{2i-1}\\
&+\frac{1}{25}L_{2i-2}+\frac{1}{25}(-1)^iL_{2i-2}-\frac{1}{50}L_{2i-1}-\frac{1}{50}(-1)^iL_{2i-1}\\
&+\frac14 F_iF_{i-1}+\frac{1}{40}F_{i-2}^2-\frac{1}{40}F_{i+1}^2-\frac{1}{50}.
\end{align*}
Applying the summation formulae
\begin{align*}
\begin{split}
&\sum_{i=2}^nL_{2i-2}L_{2i-1}=n-3+F_{4n-1}\\
&\sum_{i=2}^nL_iL_{i-1}=L_{2n}-\frac{7+(-1)^n}{2}\\
&\sum_{i=2}^niF_{2i-1}=nF_{2n+1}-(n+1)F_{2n-1}\\
&\sum_{i=2}^nF_{2i-1}=F_{2n}-1\\
&\sum_{i=2}^nL_{2i-2}=L_{2n-1}-1\\
&\sum_{i=2}^n(-1)^iL_{2i-2}=(-1)^nF_{n-1}L_n
\end{split}
\begin{split}
&\sum_{i=2}^nL_{2i-1}=L_{2n}-3\\
&\sum_{i=2}^n(-1)^iL_{2i-1}=(-1)^nF_{n-1}L_{n+1}\\
&\sum_{i=2}^nF_iF_{i-1}=F_n^2+\frac{(-1)^n-1}{2}\\
&\sum_{i=2}^nF_{i-2}^2=F_{n-1}F_{n-2}\\
&\sum_{i=2}^nF_{i+1}^2=F_{n+1}F_{n+2}-2
\end{split}
\end{align*}
leads to the identity
\begin{align*}
&\sum_{i=2}^n\sum_{j=2}^i\bigg(F_{j-1}+\sum_{k=j+1}^iF_{k-1}F_{k-j}\bigg)^2\\
&=\frac{3}{20}+\frac{(-1)^n}{8}+\frac{3}{200}(-1)^n+\frac{n}{50}+\frac{1}{40}F_{n-2}F_{n-1}+\frac{1}{4}F_n^2-\frac15 F_{2n}\\
&\quad -\frac{1}{40}F_{n+1}F_{n+2}-\frac15 F_{2n-1}-\frac15 nF_{2n-1}+\frac15 nF_{2n+1}+\frac{1}{25}F_{4n-1}\\
&\quad +\frac{1}{25}(-1)^nF_{n-1}L_n-\frac{1}{20}L_{2n}-\frac{1}{50}(-1)^nF_{n-1}L_{n+1}+\frac{1}{25}L_{2n-1}.
\end{align*}
Meanwhile, it holds that
\begin{align*}
1+\sum_{i=2}^n(1+F_iF_{i-1})^2&=n+2\sum_{i=2}^nF_iF_{i-1}+\sum_{i=2}^nF_i^2F_{i-1}^2\\
&=2F_n^2+(-1)^n-1+\frac{24}{25}n+\frac{1}{25}F_{4n}+\frac{2}{25}(-1)^nF_nL_n,
\end{align*}
since $\sum_{i=2}^nF_iF_{i-1}=F_n^2+\frac{(-1)^n-1}{2}$ and $\sum_{i=2}^nF_i^2F_{i-1}^2=-\frac{n}{25}+\frac{1}{25}F_{4n}+\frac{2}{25}(-1)^nF_nL_n$.

Putting the previous formulae together results in the equation
\begin{align}
&1+\sum_{i=2}^n(1+F_iF_{i-1})^2+2\sum_{i=2}^n\sum_{j=2}^i\bigg(F_{j-1}+\sum_{k=j+1}^iF_{k-1}F_{k-j}\bigg)^2\notag\\
&=\frac{1}{25}F_{4n}+\frac{2}{25}F_{4n-1}\label{eq:green}\\
&\quad +n-\frac25 nF_{2n-1}+\frac25 nF_{2n+1}\label{eq:red}\\
&\quad -\frac{7}{10}+\frac{2}{25}(-1)^nF_{n-1}L_n+\frac{2}{25}(-1)^nF_nL_n-\frac{1}{25}(-1)^nF_{n-1}L_{n+1}\label{eq:blue}\\
&\quad +\frac{32}{25}(-1)^n\!+\!\frac{1}{20}F_{n-2}F_{n-1}\!+\!\frac{5}{2}F_n^2\!-\!\frac25 F_{2n}\!-\!\frac{1}{20}F_{n+1}F_{n+2}\!-\!\frac25 F_{2n-1}\!-\!\frac{1}{10}L_{2n}\!+\!\frac{2}{25}L_{2n-1}.\label{eq:black}
\end{align}

At this juncture, one can proceed as follows.
\begin{itemize}
\item To simplify row~\eqref{eq:green}, use the identity
$$
\frac{1}{25}F_{4n}+\frac{2}{25}F_{4n-1}=\frac{1}{25}L_{4n}.
$$
\item To simplify row~\eqref{eq:red}, apply
$$
\frac25 nF_{2n+1}=\frac25 n(F_{2n}+F_{2n-1}).
$$
\item To cope with row~\eqref{eq:blue}, use the identity
\begin{align*}
-\frac{7}{10}+\frac{2}{25}(-1)^nF_{n-1}L_n+\frac{2}{25}(-1)^nF_nL_n-\frac{1}{25}(-1)^nF_{n-1}L_{n+1}=(-1)^n\frac{L_{2n}}{25}-\frac{33}{50}.
\end{align*}
\item Finally, to handle row~\eqref{eq:black}, utilize the identity
\begin{align*}
\frac{32}{25}(-1)^n\!+\!\frac{1}{20}F_{n-2}F_{n-1}\!+\!\frac52 F_n^2\!-\!\frac25F_{2n}\!-\!\frac{1}{20}F_{n+1}F_{n+2}\!-\!\frac25 F_{2n-1}\!-\!\frac{1}{10}L_{2n}\!+\!\frac{2}{25}L_{2n-1}\notag\\
=\frac{3}{25}L_{2n}\!+\!\frac{13(-1)^n}{50}.
\end{align*}
\end{itemize}
It is straightforward to verify the validity of each of these formulae. Altogether, the above formulae yield
\begin{align*}
&1+\sum_{i=2}^n(1+F_iF_{i-1})^2+2\sum_{i=2}^n\sum_{j=2}^i\bigg(F_{j-1}+\sum_{k=j+1}^iF_{k-1}F_{k-j}\bigg)^2\\
&=\frac{1}{25}L_{4n}+\frac{3+(-1)^n}{25}{L_{2n}}+\frac25 nF_{2n}+\frac{13(-1)^n-33}{50}+{n}.
\end{align*}
The claim follows by expanding the Fibonacci and Lucas numbers in terms of the golden ratio using~\eqref{eq:binet}.\endproof

It is evident that Theorem~\ref{thm:main} can be recast in the following way.

\begin{corollary}
For $n$ odd, it holds that
$$
c_n\geq \frac{1}{\sqrt{\frac{1}{25}\varphi^{-4n}+\frac{2}{25}\varphi^{-2n}-\frac{2}{5\sqrt{5}}n\varphi^{-2n}-\frac{23}{25}+n+\frac{2}{25}\varphi^{2n}+\frac{2}{5\sqrt{5}}n\varphi^{2n}+\frac{1}{25}\varphi^{4n}}},
$$
and for $n$ even, one has
$$
c_n\geq \frac{1}{\sqrt{\frac{1}{25}\varphi^{-4n}+\frac{4}{25}\varphi^{-2n}-\frac{2}{5\sqrt{5}}n\varphi^{-2n}-\frac{2}{5}+n+\frac{4}{25}\varphi^{2n}+\frac{2}{5\sqrt{5}}n\varphi^{2n}+\frac{1}{25}\varphi^{4n}}}.
$$
\end{corollary}

\subsection{Numerical experiments}\label{sec:numex}
The sharpness of the lower bound presented in Theorem~\ref{thm:main} is assessed by numerical experiments.  The characterization~\eqref{eq:ihm} provides an easy way of computing the numerical value of $c_n$ for $n\in\mathbb{Z}_+$. The value of the lower bound corresponding to $c_n$ is denoted by $\|Z_n\|_{\rm F}^{-1}$. 

\begin{figure}[!t]
\centering
\includegraphics[height=.313\textwidth]{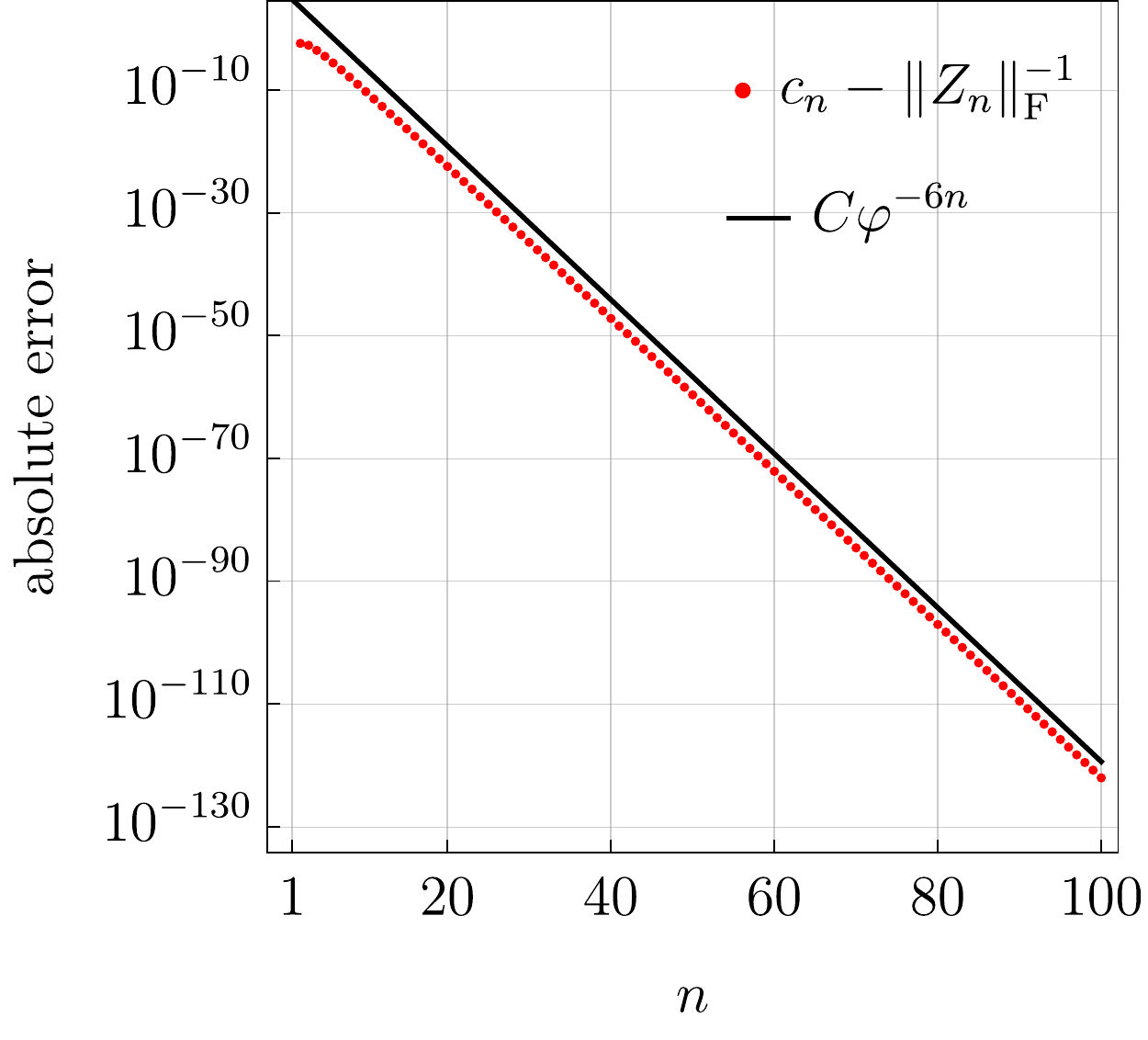}\includegraphics[height=.315\textwidth]{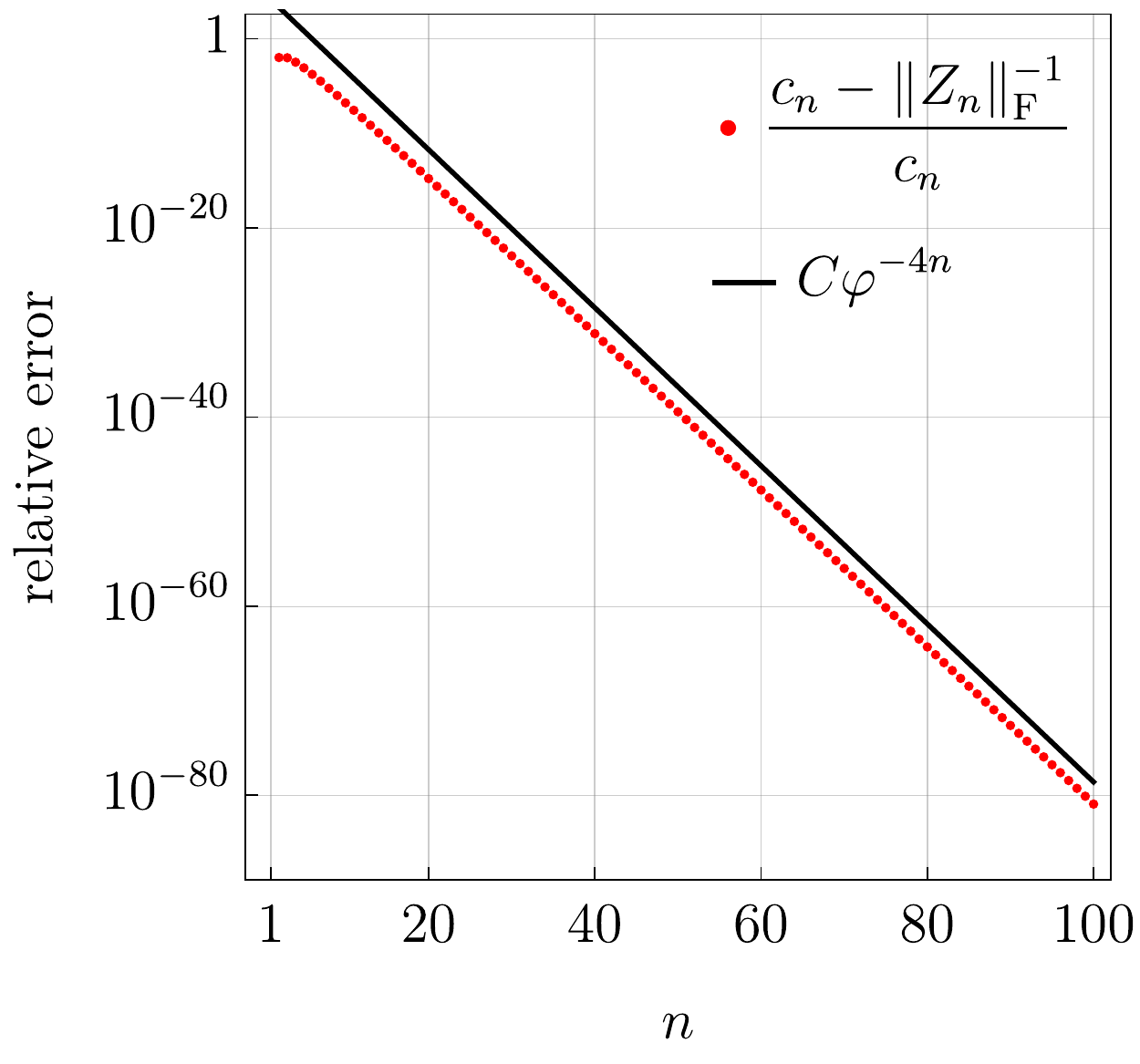}\includegraphics[height=.322\textwidth]{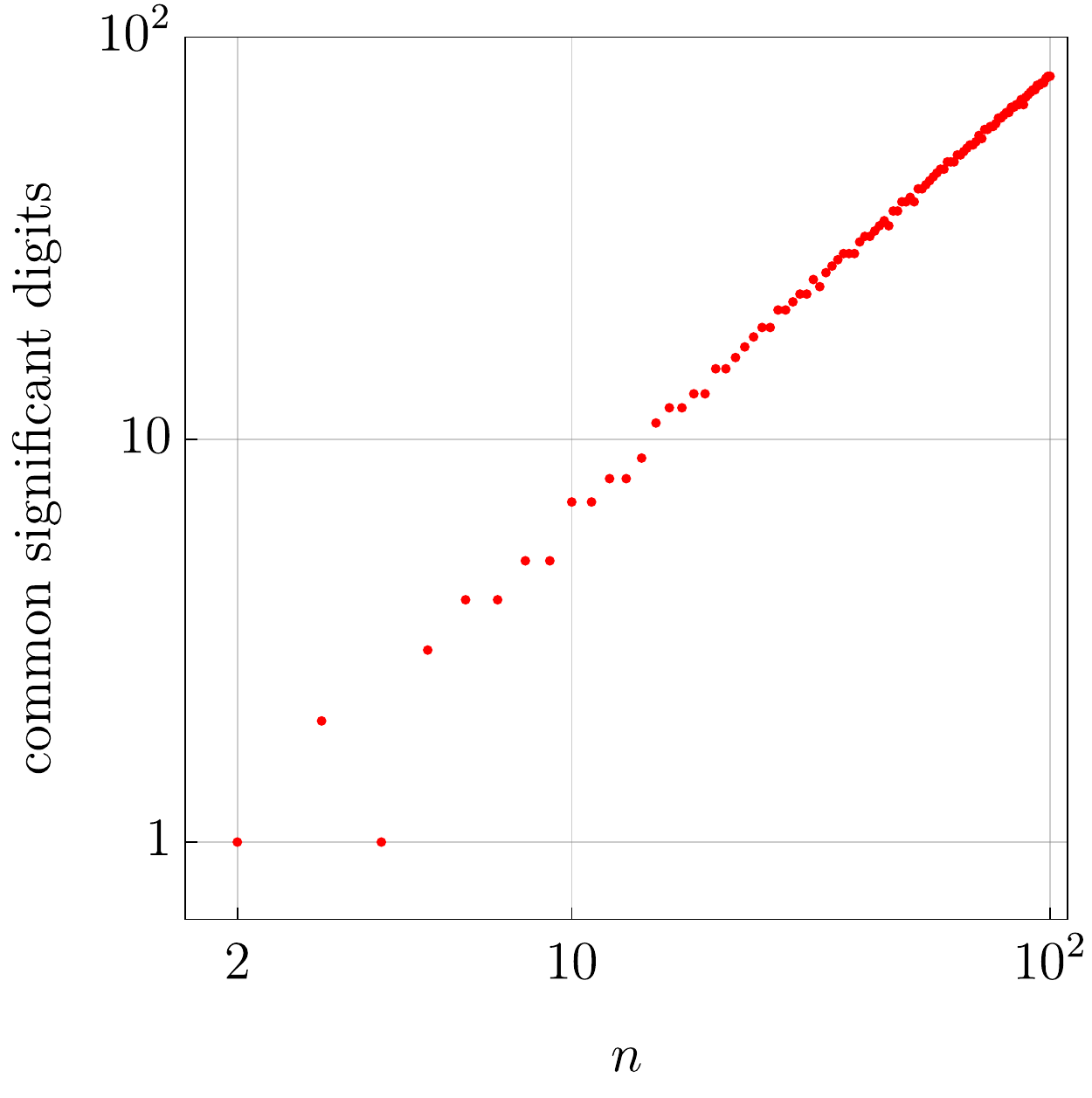}
\caption{Left and middle images: both the absolute error and the relative error between $c_n$ and the lower bound of Theorem~\ref{thm:main} decay at an exponential rate. Right image: the number of common significant digits between $c_n$ and the lower bound of Theorem~\ref{thm:main} is displayed for increasing $n$.}\label{fig:fig}
\end{figure}

In Table~\ref{table:table}, the values of both $c_n$ and the lower bound of Theorem~\ref{thm:main} have been tabulated for $n\in\{1,\ldots,10\}$. The results suggest that the lower bound becomes {\em sharper} with increasing $n$. This observation is further backed by the results illustrated in Figure~\ref{fig:fig}, where the absolute errors, relative errors as well as the number of common significant digits between $c_n$ and the lower bound in Theorem~\ref{thm:main} have been tabulated for $n\in\{2,\ldots,100\}$.

\begin{figure}[!t]
\centering
\includegraphics[height=.315\textwidth]{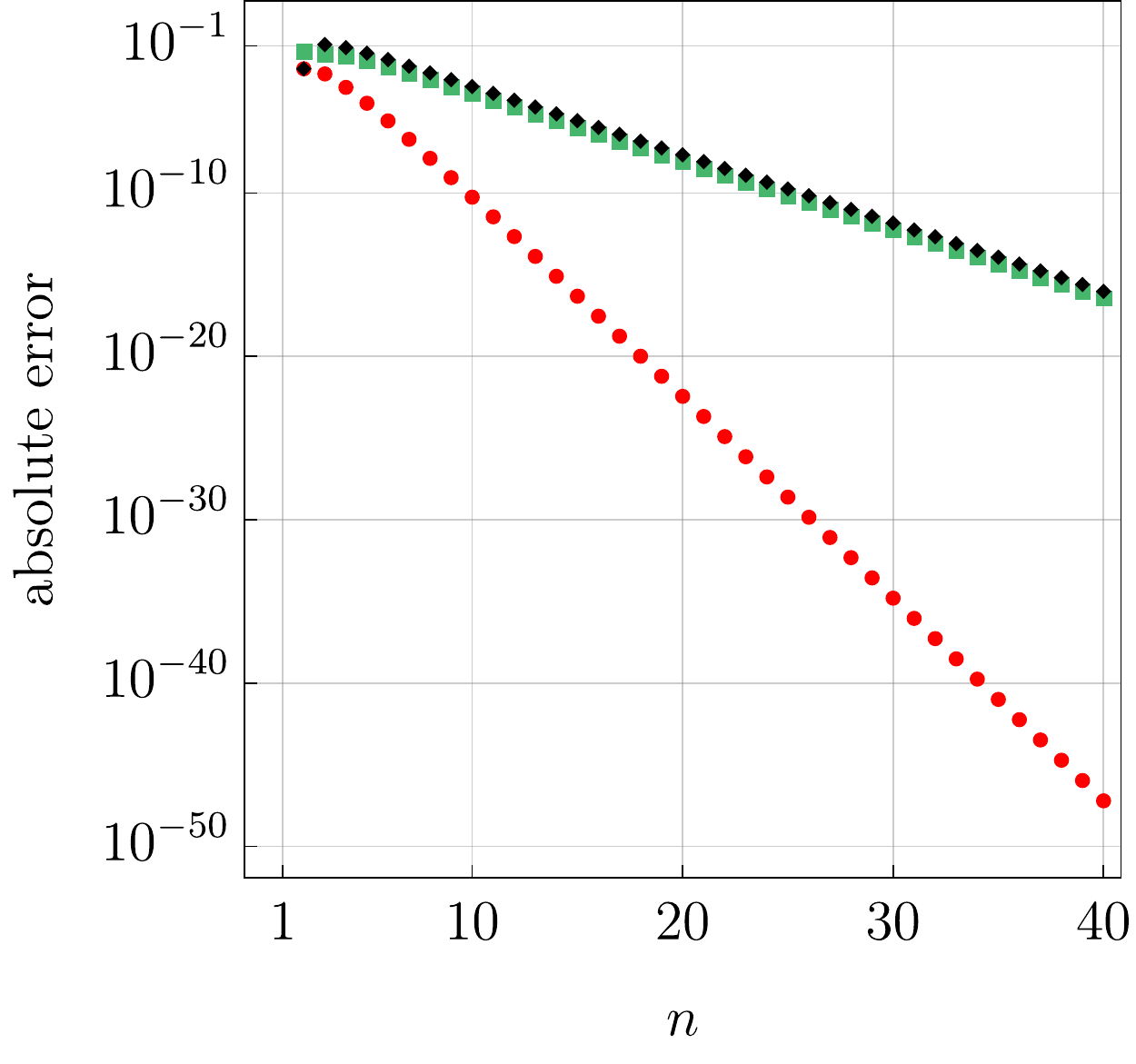}~\includegraphics[height=.315\textwidth,trim=0 11 0 11,clip]{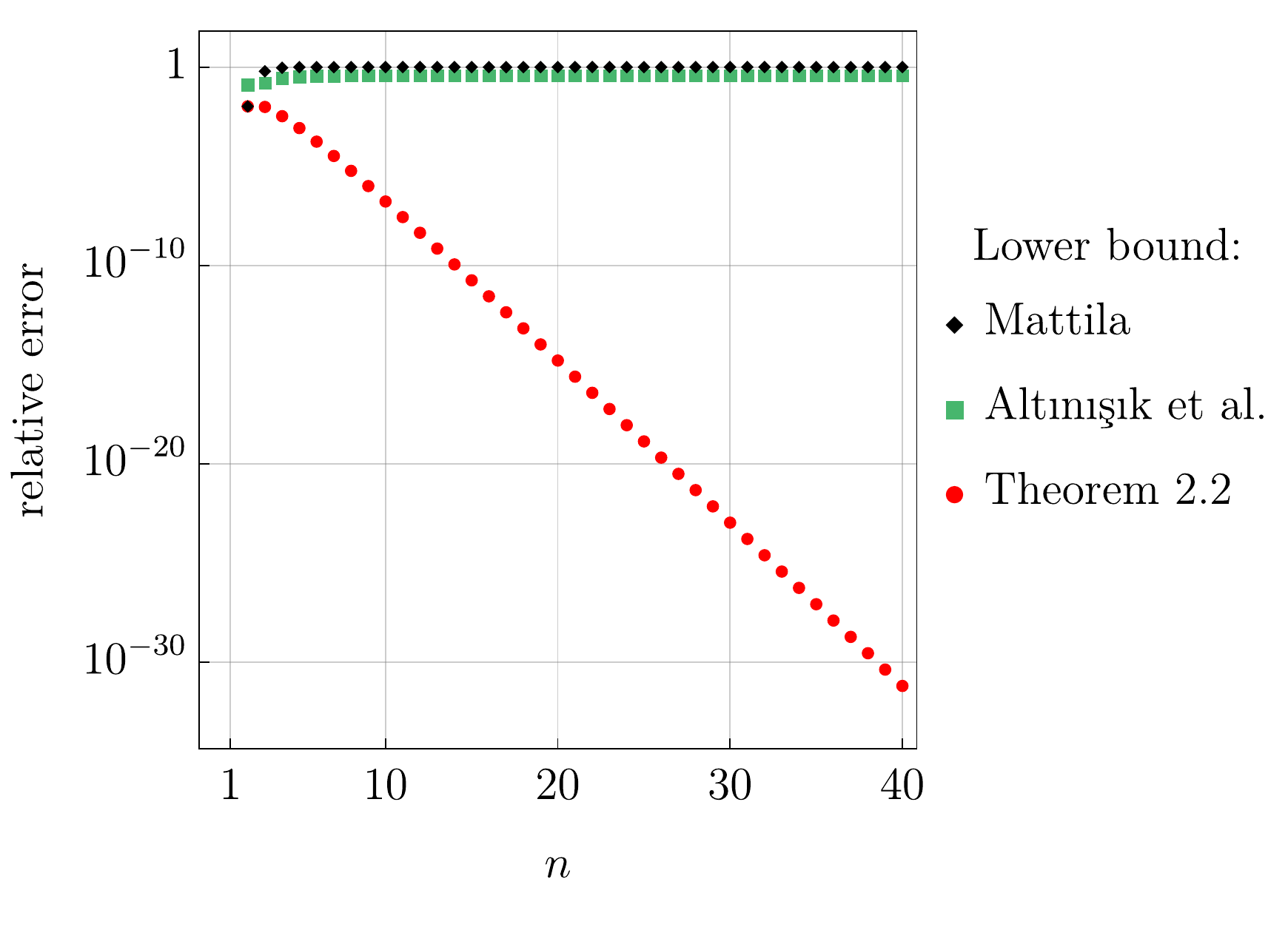}
\caption{Left image: the absolute errors between $c_n$ and Mattila's the lower bounds~\eqref{eq:mattila1}--\eqref{eq:mattila2}, $c_n$ and the lower bound \eqref{eq:altinisik} by Alt{\i}n{\i}\c{s}{\i}k et al., and $c_n$ and the lower bound of Theorem~\ref{thm:main} for increasing $n$. Right image: the corresponding relative errors.}\label{fig:fig2}
\end{figure}

The absolute and relative errors were also computed between
\begin{itemize}
\item $c_n$ and the lower bounds~\eqref{eq:mattila1}--\eqref{eq:mattila2} derived by Mattila~\cite{mattila};
\item $c_n$ and  the lower bound~\eqref{eq:altinisik} derived by Alt{\i}n{\i}\c{s}{\i}k et al.~\cite{altinisik}.
\end{itemize}
These results are displayed in Figure~\ref{fig:fig2} for $n\in\{2,\ldots,40\}$ alongside the corresponding errors between the lower bound of Theorem~\ref{thm:main} and $c_n$. The lower bounds~\eqref{eq:mattila1}--\eqref{eq:mattila2} and~\eqref{eq:altinisik} tend to zero at a faster rate than both $c_n$ and the lower bound given by~Theorem~\ref{thm:main}, explaining why the relative errors for the bounds~\eqref{eq:mattila1}--\eqref{eq:mattila2} and~\eqref{eq:altinisik} do not tend to zero as well as the poor rates of convergence in the left-hand side image of Figure~\ref{fig:fig2}.

All numerical experiments were carried out by using 150 digit precision computations in Mathematica~11.2.  

\begin{table}[!t]
\centering
\begin{tabular}{c|llll}
$n$&$c_n$&Theorem~\ref{thm:main}&Lower bound~\eqref{eq:altinisik}&Lower bounds~\eqref{eq:mattila1}--\eqref{eq:mattila2}\\
\hline $1$&$1.000000000$&$1.000000000$&$1.000000000$&$1.000000000$\\
$2$&$0.381966011$&$0.377964473$&$0.333333333$&$0.377964473$\\
$3$&$0.198062264$&$0.196116135$&$0.166666667$&$0.076923077$\\
$4$&$0.087003112$&$0.086710997$&$0.062500000$&$0.006749366$\\
$5$&$0.037068335$&$0.037037037$&$0.025000000$&$0.000540833$\\
$6$&$0.014827585$&$0.014824986$&$0.009523810$&$0.000020528$\\
$7$&$0.005816999$&$0.005816805$&$0.003663004$&$8.16298\mathrm{e}{-7}$\\
$8$&$0.002245345$&$0.002245332$&$0.001398601$&$1.62711\mathrm{e}{-8}$\\
$9$&$0.000862203$&$0.000862202$&$0.000534759$&$3.63629\mathrm{e}{-10}$\\
$10$&$0.000330004$&$0.000330004$&$0.000204248$&$4.33809\mathrm{e}{-12}$
\end{tabular}
\caption{Tabulated values of the constant $c_n$, the lower bound of Theorem~\ref{thm:main}, the lower bound~\eqref{eq:altinisik} derived by Alt{\i}n{\i}\c{s}{\i}k et al.~\cite{altinisik}, and the lower bounds~\eqref{eq:mattila1}--\eqref{eq:mattila2} derived by Mattila~\cite{mattila} for $n\in\{1,\ldots,10\}$. The numerical values suggest that the lower bound of Theorem~\ref{thm:main} becomes sharper as $n$ increases, while the other bounds become less sharp with increasing $n$. Note that the lower bound given by Theorem~\ref{thm:main} coincides with Mattila's bound for $n\in\{1,2\}$.}\label{table:table}
\end{table}


\section{The Ilmonen--Haukkanen--Merikoski numbers $C_n$}\label{sec:ihm}

To conclude this paper, the following new characterization is proved for the Ilmonen--Haukkanen--Merikoski numbers $C_n$.
\begin{theorem}\label{lemma:ihm} It holds for all $n\in\mathbb{Z}_+$ that
\begin{align*}
C_n=\frac14\csc^2\bigg(\frac{\pi}{4n+2}\bigg)=\frac{4n^2}{\pi^2}+\frac{4n}{\pi^2}+\bigg(\frac{1}{12}+\frac{1}{\pi^2}\bigg)+\mathcal{O}\bigg(\frac{1}{n^2}\bigg).
\end{align*}
\end{theorem}
\proof

The second identity is a consequence of the Laurent expansion of the first expression developed at infinity. It is therefore enough to focus on proving the first identity.

It is easy to check that the claim holds for $n=1$. Let $n\geq 2$. By~\cite{ihm}, it is known that
$$
C_n=\lambda_{\max}(W_n),
$$
where $W_n$ is the $n\times n$ matrix defined elementwise by
$$
(W_n)_{i,j}=\min\{i,j\},\quad i,j\in\{1,\ldots,n\}.
$$
It is easy to see that its matrix inverse $B_n=W_n^{-1}$ is the tridiagonal matrix
$$
\quad B_n=\begin{pmatrix}
2&-1&&&&\\
-1&2&-1&&&\\
&-1&2&-1&&\\
&&&\ddots&&\\
&&&-1&2&-1\\
&&&&-1&1
\end{pmatrix}
$$
and hence $C_n=\lambda_{\min}(B_n)^{-1}$. In other words, it is sufficient to find the reciprocal of the minimal eigenvalue of $B_n$. Noting that $B_n$ is a special case of a second order finite difference matrix subject to mixed Dirichlet--Neumann boundary conditions (see also the Remark following this proof), it is well known that the eigenvalues of $B_n$ are roots of certain Chebyshev polynomials and, as such, the roots have closed form solutions. A brief derivation is presented in the following for completeness.

Let $A_n=B_n+{e}_n{e}_n^{\rm T}$, where ${e}_n=[0,0,\ldots,0,1]^{\rm T}\in \mathbb{R}^n$. Let $p_n(\lambda)=\det(A_n-\lambda I_n)$ and $q_n(\lambda)=\det(B_n-\lambda I_n)$ be the characteristic polynomials of $A_n$ and $B_n$, respectively. The matrix $A_n$ is a tridiagonal Toeplitz matrix and it has the eigenvalues (cf., e.g.,~\cite{triagonaltoeplitz})
$$
\mu_k=2-2\cos\bigg(\frac{\pi k}{n+1}\bigg),\quad k\in \{1,\ldots,n\}.
$$

It is useful to consider {\em Chebyshev polynomials of the second kind} $U_n$, which are characterized by the three-term recurrence
\begin{align}
U_0(x)= 1,\quad U_1(x)=2x,\quad \text{and}\quad U_n(x)=2xU_{n-1}(x)-U_{n-2}(x)\quad\text{for }n\geq 2.\label{eq:threeterm}
\end{align}
If $|x|<1$, then it additionally holds that (cf., e.g.,~\cite{gautschi})
\begin{align}
U_n(x)=\frac{\sin((n+1)\,{\rm arc\,cos}(x))}{\sqrt{1-x^2}}.\label{eq:arcsin}
\end{align}
Equation~\eqref{eq:arcsin} implies that the roots of $U_n$ lie in the open interval $(-1,1)$ and they are given explicitly by the formula
$$
x_k=\cos\bigg(\frac{\pi k}{n+1}\bigg),\quad k\in \{1,\ldots,n\}.
$$
Hence $\mu_k=2-2x_k$ for all $k\in\{1,\ldots,n\}$. Since $p_n$ is a monic polynomial, it follows that
$$
p_n(\lambda)=\prod_{k=1}^n(\mu_k-\lambda)=\prod_{k=1}^n(2-2x_k-\lambda)=2^n\prod_{k=1}^n\bigg(\bigg(1-\frac{\lambda}{2}\bigg)-x_k\bigg)=U_n\bigg(1-\frac{\lambda}{2}\bigg),
$$
where the final equality follows from the identity $U_n(x)=2^n\prod_{k=1}^n(x-x_k)$ derived using the fundamental theorem of algebra, the leading coefficient being a consequence of the three-term recurrence~\eqref{eq:threeterm}.

Developing the Laplace cofactor expansion of $\det(B_n-\lambda I_n)$ across the final column and using the properties~\eqref{eq:threeterm} and~\eqref{eq:arcsin} yields that
\begin{align*}
q_n(\lambda)&=(1-\lambda)p_{n-1}(\lambda)-p_{n-2}(\lambda)\\
&=2\bigg(1-\frac{\lambda}{2}\bigg)U_{n-1}\bigg(1-\frac{\lambda}{2}\bigg)-U_{n-2}\bigg(1-\frac{\lambda}{2}\bigg)-U_{n-1}\bigg(1-\frac{\lambda}{2}\bigg)\\
&=U_{n}\bigg(1-\frac{\lambda}{2}\bigg)-U_{n-1}\bigg(1-\frac{\lambda}{2}\bigg)\\
&=\frac{\sin\big((n+1)\, {\rm arc\,cos}\big(1-\frac{\lambda}{2}\big)\big)-\sin\big(n\, {\rm arc\,cos}\big(1-\frac{\lambda}{2}\big)\big)}{\sqrt{1-(1-\frac{\lambda}{2})^2}}.
\end{align*}
The previous expression can be used to solve the roots of $q_n$ by elementary means, i.e.,
$$
\lambda_j=4\cos^2\bigg(\frac{j\pi}{2n+1}\bigg),\quad j\in\{1,\ldots,n\}.
$$
Since the smallest root of $q_n$ is $\lambda_n$ for all $n\in\mathbb{Z}_+$, it follows that $$C_n=\frac{1}{\lambda_{\min}(B_n)}=\frac14 \sec^2\bigg(\frac{n\pi}{2n+1}\bigg)=\frac14\csc^2\bigg(\frac{\pi}{4n+2}\bigg)$$
completing the proof.\endproof

\emph{Remark.} The matrix $B_n$ is (up to a scalar multiple) precisely the finite difference matrix corresponding to the Dirichlet--Neumann problem
$$
-u''(x)=f(x)\quad \text{for }x\in (a,b),\quad u(a)=0,\quad u'(b)=0.
$$
The properties of finite difference matrices for this problem are very well known in the literature; see, e.g.,~\cite{fde} for a comprehensive treatment of the topic. 

It was shown in~\cite{ihm} that the numbers $C_n$ can be bounded by
$$
C_n\leq \sqrt{(2n-1)+4(2n-3)+9(2n-5)+\cdots+3(n-1)^2+n^2},\quad n\in\mathbb{Z}_+,
$$
but the closed form solution stated in Theorem~\ref{lemma:ihm} appears to have eluded the authors of the aforementioned paper.

\section*{Conclusions}
The numerical experiments presented in this paper suggest that the lower bound obtained for the numbers $c_n$ is extremely sharp as $n$ tends to infinity. Theorems~\ref{thm:main} and~\ref{lemma:ihm} can be used in conjunction with eigenvalue bounds such as~\eqref{eq:application} to obtain new explicit lower and upper bounds for the smallest and largest eigenvalues of power GCD matrices as well as lattice-theoretic meet and join matrices. See~\cite[Sections~3--6]{ihm} for further discussion on eigenvalue bounds involving $c_n$~and $C_n$ as well as~\cite[Theorem~4.2]{hongloewy} for the special case of power GCD matrices.

The numerical evidence leads the author to further conjecture that $c_n\sim 5\varphi^{-2n}$ as $n\to\infty$, based on the dominating term that appears in Theorem~\ref{thm:main}. Proving this asymptotic result appears to require developing mathematical techniques which are beyond the scope of this paper, posing an interesting challenge for researchers working in this area.

\section*{Acknowledgements}
The author gratefully acknowledges the financial support from the Australian Research Council (grant number DP180101356). The author thanks the anonymous referees for their valuable comments and suggestions which helped to improve this paper.

\bibliographystyle{plain}
\bibliography{boolean}

\begin{thebibliography}{10}

\bibitem{altinisik15}
E.~Alt{\i}n{\i}\c{s}{\i}k and \c{S}. B\"{u}y\"{u}kk\"{o}se.
\newblock A proof of a conjecture on monotonic behavior of the smallest and the
  largest eigenvalues of a number theoretic matrix.
\newblock {\em Linear Algebra Appl.}, 471:141--149, 2015.

\bibitem{altinisik16}
E.~Alt{\i}n{\i}\c{s}{\i}k and \c{S}. B\"{u}y\"{u}kk\"{o}se.
\newblock On bounds for the smallest and the largest eigenvalues of {GCD} and
  {LCM} matrices.
\newblock {\em Math. Inequal. Appl.}, 19(1):117--125, 2016.

\bibitem{altinisik}
E.~Alt{\i}n{\i}\c{s}{\i}k, A.~Keskin, M.~Y{\i}ld{\i}z, and M.~Demirb\"{u}ken.
\newblock On a conjecture of {I}lmonen, {H}aukkanen and {M}erikoski concerning
  the smallest eigenvalues of certain {G}{C}{D} related matrices.
\newblock {\em Linear Algebra Appl.}, 493:1--13, 2016.

\bibitem{balatoni}
F.~Balatoni.
\newblock On the eigenvalues of the matrix of the {S}mith determinant.
\newblock {\em Mat. Lapok}, 20:397--403, 1969.
\newblock In Hungarian.

\bibitem{beslinligh}
S.~Beslin and S.~Ligh.
\newblock Greatest common divisor matrices.
\newblock {\em Linear Algebra Appl.}, 118:69--76, 1989.

\bibitem{bourqueligh}
K.~Bourque and S.~Ligh.
\newblock Matrices associated with arithmetical functions.
\newblock {\em Linear Multilinear Algebra}, 34:261--267, 1993.

\bibitem{gautschi}
W.~Gautschi.
\newblock {\em Orthogonal Polynomials: Computation and Approximation}.
\newblock Oxford University Press, 2004.

\bibitem{inertia}
P.~Haukkanen, M.~Mattila, and J.~M\"{a}ntysalo.
\newblock Studying the inertias of {LCM} matrices and revisiting the
  {B}ourque-{L}igh conjecture.
\newblock {\em J. Combin. Theory Ser. A}, 171:105161, 2020.

\bibitem{hong2008}
S.~Hong.
\newblock Asymptotic behavior of largest eigenvalue of matrices associated with
  completely even functions (mod $r$).
\newblock {\em Asian-Eur. J. Math.}, 1(2):225--235, 2008.

\bibitem{hongenochlee2008}
S.~Hong and K.~S. {Enoch Lee}.
\newblock Asymptotic behavior of eigenvalues of reciprocal power {LCM}
  matrices.
\newblock {\em Glasgow Math. J.}, 50:163--174, 2008.

\bibitem{hongloewy}
S.~Hong and R.~Loewy.
\newblock Asymptotic behavior of eigenvalues of greatest common divisor
  matrices.
\newblock {\em Glasgow Math. J.}, 46(3):551--569, 2004.

\bibitem{hongloewy11}
S.~Hong and R.~Loewy.
\newblock Asymptotic behavior of the smallest eigenvalue of matrices associated
  with completely even functions (mod $r$).
\newblock {\em Int. J. Number Theory}, 7(6):1681--1704, 2011.

\bibitem{ilmonen}
P.~Ilmonen.
\newblock On meet hypermatrices and their eigenvalues.
\newblock {\em Linear Multilinear Algebra}, 64(5):842--855, 2016.

\bibitem{ihm}
P.~Ilmonen, P.~Haukkanen, and J.~K. Merikoski.
\newblock On eigenvalues of meet and join matrices associated with incidence
  functions.
\newblock {\em Linear Algebra Appl.}, 429:859--874, 2008.

\bibitem{mattila}
M.~Mattila.
\newblock On the eigenvalues of combined meet and join matrices.
\newblock {\em Linear Algebra Appl.}, 466:1--20, 2015.

\bibitem{haukkanenmattila12}
M.~Mattila and P.~Haukkanen.
\newblock On the eigenvalues of certain number-theoretic matrices.
\newblock {\em East-West J. Math.}, 14(2):121--130, 2012.

\bibitem{merikoski}
J.~K. Merikoski.
\newblock Lower bounds for the largest eigenvalue of the {GCD} matrix on
  $\{1,2,\ldots,n\}$.
\newblock {\em Czechoslovak Math. J.}, 66:1027--1038, 2016.

\bibitem{fde}
A.~R. Mitchell and D.~F. Griffiths.
\newblock {\em The Finite Difference Method in Partial Differential Equations}.
\newblock John Wiley \& Sons, 1980.

\bibitem{triagonaltoeplitz}
S.~Noschese, L.~Pasquini, and L.~Reichel.
\newblock Tridiagonal {T}oeplitz matrices: properties and novel applications.
\newblock {\em Numer. Linear Algebra Appl.}, 20(2):302--326, 2013.

\bibitem{smith}
H.~J.~S. Smith.
\newblock On the value of a certain arithmetical determinant.
\newblock {\em Proc. London Math. Soc.}, 7:208--212, 1875/76.

\end{thebibliography}
\end{document}